\documentclass[reqno,11pt]{amsart}

\usepackage{booktabs}
\usepackage{cases}
\usepackage{mathrsfs}
\usepackage[T1]{fontenc}
\usepackage{mathrsfs}
\usepackage{amsmath,latexsym,amssymb,amsfonts,amsbsy, amsthm}
\usepackage{bm}
\usepackage[usenames]{color}
\usepackage{xspace,colortbl}
\usepackage{epsfig}
\usepackage{graphicx}
\usepackage{subfigure}
\usepackage{amsmath,amsfonts,amsthm,amssymb,amscd}
\usepackage{color}
\usepackage[title,titletoc,toc]{appendix}
\usepackage{bigints}
\usepackage{hyperref}
\usepackage{cases}\usepackage{extarrows}
\usepackage{setspace}
\usepackage{titletoc}
\usepackage{color}
\usepackage{caption}
\usepackage[numbers,sort&compress]{natbib}
\usepackage{amsmath}
\usepackage{graphicx}
\usepackage{tikz}
\usepackage{xcolor}
\usepackage{pgfplots}
\pgfplotsset{compat=1.18}

\usepackage{tikz}
\usepackage{pgfplots}
\pgfplotsset{compat=1.18}
\usepgfplotslibrary{colormaps}
\usepackage{amsmath}
\usepackage{amssymb}
\usepackage{tikz}
\usepackage{pgfplots}
\pgfplotsset{compat=1.18}
\usepackage{amsmath,amssymb}
\usepackage{lipsum} 
\usepackage{tikz}
\usetikzlibrary{arrows.meta, shapes.geometric, positioning}
\usepackage{tikz}
\usetikzlibrary{decorations.pathreplacing}
\usetikzlibrary{shapes.geometric, decorations.pathreplacing, shadows}
\usepackage{tikz-3dplot}
\usetikzlibrary{shapes.geometric, calc, decorations.pathmorphing}

\usepgfplotslibrary{patchplots}

\usepackage{tikz}
\usepackage{amsmath, amssymb}
\usetikzlibrary{shadows, calc, arrows.meta, positioning, shapes.arrows}
\definecolor{topcolor}{RGB}{230,230,230}
\definecolor{bottomcolor}{RGB}{150,150,150}
\usepackage{pgfplots}
\pgfplotsset{compat=1.18}
\usepackage{pgfplots}
\pgfplotsset{compat=1.18}

\usepackage{rotating} 
\usepackage[most]{tcolorbox}

\usetikzlibrary{arrows.meta}
\usetikzlibrary{shapes,arrows.meta,positioning}
\usetikzlibrary{arrows.meta, positioning}

\usepackage{amsmath}
\usepackage{amsfonts}

\usepackage{amsmath}
\usepackage[most]{tcolorbox}

\usepackage[most]{tcolorbox}

\tcbset{
    mybox/.style={
        colframe=black, 
        colback=#1,     
        boxrule=0.4pt,  
        arc=0pt,        
        left=0pt, right=0pt, top=0pt, bottom=0pt, 
        boxsep=0pt,     
        fontupper=\footnotesize, 
        nobeforeafter,  
        hbox           
    }
}

\tikzstyle{block} = [rectangle, draw, fill=blue!20, 
    text width=5em, text centered, rounded corners, minimum height=4em, font=\small]
\tikzstyle{description} = [rectangle, draw=none, fill=yellow!20, 
    text width=8em, text centered, minimum height=3em, font=\footnotesize, inner sep=2pt]
\tikzstyle{line} = [draw, -{Latex[length=2mm]}, thick]

\tikzset{
    block/.style = {rectangle, draw, text width=6em, text centered, rounded corners, minimum height=4em},
    line/.style = {draw, -Latex},
    description/.style = {text width=8cm, align=left, font=\small},
    column sep=4cm, 
    row sep=1.5cm  
}

\definecolor{mainColor}{RGB}{213, 94, 0} 
\definecolor{secondColor}{RGB}{0, 158, 115} 
\definecolor{thirdColor}{RGB}{86, 180, 233} 
\definecolor{fourthColor}{RGB}{230, 159, 0} 

\definecolor{mainColor}{HTML}{0073CF} 
\definecolor{secColor}{HTML}{1F4E79} 

\definecolor{mainColor}{HTML}{4CAF50} 
\definecolor{secColor}{HTML}{2196F3}  
\definecolor{terColor}{HTML}{FF9800}  
\definecolor{quaColor}{HTML}{F44336}  

\usetikzlibrary{shapes,arrows,positioning,backgrounds}
\pgfplotsset{compat=1.18} 

\usetikzlibrary{shapes, arrows.meta, positioning, fit, backgrounds}

\tikzset{
    block/.style = {rectangle, draw, text width=8em, text centered, rounded corners, minimum height=4em, font=\small, align=center, fill=blue!20},
    detail/.style = {rectangle, draw, text width=8em, text centered, rounded corners, minimum height=4em, font=\footnotesize, align=center, fill=green!20},
    line/.style = {draw, -Latex, thick, blue},
    connection/.style = {draw, -Latex, thick, red},
    cloud/.style = {ellipse, draw, text width=6em, text centered, minimum height=3em, font=\small, align=center, fill=yellow!20, drop shadow={shadow xshift=0.5ex,shadow yshift=-0.5ex,fill=gray,opacity=0.2}},
    title/.style = {font=\bfseries, align=center, font=\large},
}

\allowdisplaybreaks

\numberwithin{equation}{section}

\newtheorem{lemma}{Lemma}[section]

\newtheorem{theorem}{Theorem}[section]

\begin{document}


\title{On a Conjecture by Ren and Li}

\author{
  Yongbing Luo and Ping Yan}
\maketitle

\newcommand{\Addresses}{{
  \small

  {\sc Yongbing Luo }\\ \textsc{College of Mathematics and Computer Science, Zhejiang A\&F University, Hangzhou, 311300, People’s Republic of China}

  \par\nopagebreak
  \textit{E-mail address:}  \texttt{luoyongbingah@163.com; luoyongbing@zafu.edu.cn}

}

{\sc Ping Yan (Corresponding author)}\\ \textsc{College of Mathematics and Computer Science, Zhejiang A\&F University, Hangzhou, 311300, People’s Republic of China}

\noindent  
{\sc Department of Mathematics and Statistics,
University of Helsinki, FIN-00014 Helsinki, Finland}

  \par\nopagebreak
  \textit{E-mail address:}  \texttt{ping.yan@helsinki.fi}

}

\begin{abstract}
 This paper proves a conjecture proposed by Ren and Li (2015: 393, \emph{Journal of Inequalities and Applications}). 
Our result eliminates the constraints on the parity and size of $m$, as well as the restriction $x > 1$, 
required in Ren and Li's theorem. Consequently, it fully subsumes their results while extending validity 
to all integers $m \geq 1$ and all $x > 0$. Crucially, we establish the inequality 
$S_m(x) > \sigma_m(x)$  unconditionally, requiring no parity conditions, 
size conditions on $m$, or lower bound on $x$.

\smallskip
\noindent {\sc 2020 Mathematics Subject Classification}: { 42B25;
42B35.}

\smallskip
\noindent {\sc Keywords\/}: {  Carleman inequeality;
integral; series; number e.
}
\end{abstract}
\section{Introduction}

\setcounter{equation}{0}

The classical Carleman inequality \cite{L-1},
$$
\sum_{n=1}^{\infty} (a_1 a_2 \cdots a_n)^{1/n} < e \sum_{n=1}^{\infty} a_n,
$$
where \( a_n \geq 0 \) for \( n = 1, 2, \dots \) and \( 0 < \sum_{n=1}^{\infty} a_n < \infty \), has been the subject of extensive study in recent years \cite{L-5,L-3,L-4,L-6}. Yang generalized a strengthened version of Carleman's inequality obtained in \cite{L-5} and conjectured that if
\begin{equation}\label{LLL-1}
\left(1 + \frac{1}{x}\right)^x = e \left(1 - \sum_{k=1}^{\infty} \frac{b_k}{(x+1)^k}\right), \ x > 0,
\end{equation}
then $ b_k > 0 $ for $k = 1, 2, \dots$. This conjecture was subsequently proven by  Gyllenberg and Yan \cite{L-3}. As an application, they established that for any positive integer $m$,
\begin{equation}\label{LLL-2}
\sum_{n=1}^{\infty} (a_1 a_2 \cdots a_n)^{1/n} < e \sum_{n=1}^{\infty} \left(1 - \sum_{k=1}^{m} \frac{b_k}{(n+1)^k}\right) a_n
\end{equation}
holds, where $a_n \geq 0 $ for $ n = 1, 2, \dots $, $ 0 < \sum\limits_{n=1}^{\infty} a_n < \infty$, and the coefficients $ b_k$ are defined recursively by
\begin{align*}
b_1 = \frac{1}{2}, 
b_{n+1} = \frac{1}{n+1} \left( \frac{1}{n+2} - \sum_{k=1}^{n} \frac{b_k}{n+2-k} \right).
\end{align*} In the concluding remarks of \cite{L-5}, Yang suggested that improved results might be obtained by replacing the right-hand side of \eqref{LLL-1} with $ e \left(1 - \sum\limits_{n=1}^{\infty} \frac{d_n}{(x+\varepsilon)^n} \right)$, where $ \varepsilon \in (0,1] $ and \( d_n = d_n(\varepsilon) \), although specific values for  $\varepsilon $ were not provided. Building upon this, Hu and Mortici \cite{L-4} demonstrated that $\varepsilon = \frac{11}{12}$ yields a series
$
\sum_{n=1}^{\infty} \frac{d_n}{(x + \tfrac{11}{12})^n}
$
with faster convergence. Consequently, they derived the inequality
\begin{equation}\label{HM-ineq}
\sum_{n=1}^{\infty} (a_1 a_2 \cdots a_n)^{1/n} < e \sum_{n=1}^{\infty} \left(1 - \sum_{k=1}^{m} \frac{d_k}{(n + \tfrac{11}{12})^k}\right) a_n,
\end{equation}
which refines \eqref{LLL-2}. Building upon this,   Ren and Li \cite{L-6} proved the following result:

\begin{theorem}[\cite{L-6}]\label{444}
Let
$$
\left(1 + \frac{1}{x}\right)^x = e \left(1 - \sum_{k=1}^{\infty} \frac{b_k}{(1+x)^k}\right) = e \left(1 - \sum_{k=1}^{\infty} \frac{d_k}{(x + \tfrac{11}{12})^k}\right), \quad x > 0,
$$
and define the partial sums
\begin{align}\label{LLLL-1}
\sigma_m(x) := \sum_{k=1}^{m} \frac{b_k}{(1+x)^k}\ \text{and}\
S_m(x) := \sum_{k=1}^{m} \frac{d_k}{(x + \tfrac{11}{12})^{k}}.
\end{align}
Then:
\begin{enumerate}
    \item[(i)] If \( m \geq 6 \) is even, \( S_m(x) > \sigma_m(x) \) for all \( x > 0 \).
    \item[(ii)] If \( m \geq 7 \) is odd, \( S_m(x) > \sigma_m(x) \) for all \( x > 1 \).
\end{enumerate}
\end{theorem}
Ren and Li also formulated the conjecture that
\begin{equation}\label{RL-conj}
S_m(x) > \sigma_m(x)\ \text{for all}\ x > 0\ \text{and all integers}\ m \geq 1.
\end{equation}
In this paper, we prove this conjecture in the following theorem.
\begin{theorem}\label{thm1} Suppose that all conditions in Theorem \ref{444} are satisfied.
Then the ineqiality \ref{RL-conj} holds.
\end{theorem}
Compared with the theorem in Ren and Li \cite{L-6}, Theorem \ref{thm1} in this paper achieves at least the following two improvements in terms of both conditions and conclusions.\\
{\bf Improvements in Conditions (Weaker Conditions in Theorem \ref{thm1})}
\begin{itemize}
    \item Conditions in  Theorem \ref{444}:
    \begin{itemize}
        \item Restrictions on  parity and size of $m$: Requires $m \geq 6$ (even) or $m \geq 7$ (odd). Restrictions on  range of $x$: For odd $m \geq 7$, requires $x > 1$ (while for even $m$, $x > 0$). These conditions limit the applicability of Theorem \ref{444}, for example:
    When $m = 1,2,3,4,5$ (all integers) or $m = 6$ (odd), the theorem does not apply. When $m$ is odd and $0 < x \leq 1$ (e.g., $x = 0.5$), the theorem does not apply.  
    \end{itemize}   
\end{itemize}
\noindent
{\bf Improvements in Conclusions (Stronger Conclusions in Theorem \ref{thm1})}
\begin{itemize}
    \item Stronger conclusions in Theorem \ref{thm1}
    \begin{itemize}
        \item Conclusions fully cover all $m \geq 1$ and all $x > 0$: including regions not covered by Theorem \ref{444} (e.g., $m = 1,2,3,4,5$; odd $m$ with $0 < x \leq 1$). The stronger conclusions are reflected in two key extensions: specifically, the range of $m$ is extended from requiring $m \geq 6$ (even) or $m \geq 7$ (odd) to covering all $m \geq 1$, and for odd $m$, the range of $x$ is extended from requiring $x > 1$ to covering all $x > 0$, especially including the case $x \leq 1$. For example, when $m = 5$ (odd, less than 7), Theorem \ref{444} does not apply, but Theorem \ref{thm1} establishes $S_5(x) > \sigma_5(x)$ for all $x > 0$; similarly, when $m = 7$ (odd) and $x = 0.5$ (less than 1), Theorem \ref{444} does not apply, but Theorem \ref{thm1} still holds.
    \end{itemize}
\end{itemize}

\section{Lemmas}
To prove our main results, we require the following lemma. Throughout this paper, we define the function \( g: [0, 1] \to \mathbb{R} \) by
\begin{equation}
g(s) := 
\begin{cases} 
\dfrac{1}{\pi} s^{s} (1-s)^{1-s} \sin(\pi s) & \text{if } 0 < s < 1 \\ 
0 & \text{if } s = 0 \text{ or } s = 1 
\end{cases}
\end{equation}
\begin{lemma}[\cite{L-6}]\label{lemma1}
For $x>0$, let
\begin{align*}
(1+x)^{x}=e\left(1-\displaystyle\sum^{\infty}_{k=1}\frac{b_k}{(1+x)^{k}}\right).
\end{align*}
Then
\begin{align}\label{0620-1}
\begin{cases}
b_{k}>0,\ k=1,2,\cdots,\\
b_{1}=\frac{1}{2},\\
b_{n+1}=\frac{1}{n+1}\left(\frac{1}{n+2}-\sum^{n}_{k=1}\frac{b_{j}}{n+2-j}\right),\ n=1,2,\cdots,\\
eb_{k}=\int^{1}_{0}g(s)s^{k-2}{\rm d}s,\ k=2,3,\cdots.
\end{cases}
\end{align}
\noindent {\bf Remark 1}. By \eqref{0620-1} (iv), we have
\begin{align*}
\begin{cases}
\int^{1}_{0}g(s)s^{n-2}{\rm d}s=\int^{1}_{0}g(s)(1-s)^{n-2}{\rm d}s=eb_{n} (n=2,3,\cdots),\\
\int^{1}_{0}g(s){\rm d}s=eb_{2}=\frac{e}{24},\  \int^{1}_{0}g(s)s{\rm d}s=eb_{3}=\frac{e}{48},\\
\int^{1}_{0}\frac{1}{s}g(s){\rm d}s=\int^{1}_{0}\frac{1}{1-s}g(s){\rm d}s
                              =\int^{1}_{0}(1+s+s^{2}+\cdots)g(s){\rm d}s\\
                           \ \ \ \  \ \ \   \ \ \   \   \ \ \     =e\displaystyle\sum^{\infty}_{n=2}b_{n}=e\displaystyle\sum^{\infty}_{n=1}b_{n}-eb_{1}
                              =e(1-\frac{1}{e})-\frac{e}{2}=\frac{e}{2}-1.
\end{cases}
\end{align*}
\end{lemma}

\begin{lemma}[\cite{L-6}]\label{lemma2}
For $x>0$, let
\begin{equation*}
(1+\frac{1}{x})^{x}=e\left(1-\displaystyle\sum^{\infty}_{k=1}\frac{d_k}{(x+\frac{11}{12})^{k}}\right).
\end{equation*}
Then
\begin{equation*}
d_{1}=\frac{1}{2}
\end{equation*}
\begin{equation*}
d_{n}=\frac{1}{12^{n-1}e}\left((-1)^{n-1}+\int^{1}_{0}\frac{(12t-1)^{n-1}}{t}g(t){\rm d}t\right),~~n=2,3,\cdots.
\end{equation*}
\end{lemma}

\begin{lemma}[\cite{L-6}]\label{lemma3}
Let
\begin{equation*}
h(x)=(1+x)\left(e-(1+\frac{1}{x})^{x}\right), x>0,
\end{equation*}
Then we have
\begin{equation*}
h(x)=\frac{e}{2}+\frac{1}{\pi}\int^{1}_{0}\frac{s^{s}(1-s)^{1-s}sin(\pi s)}{x+s}{\rm d}s.
\end{equation*}
\end{lemma}

\section{Proof of Theorem \ref{thm1}}
{\bf Proof of Theorem \ref{thm1}.} By (\ref{LLLL-1}), Lemma \ref{lemma1} and Lemma \ref{lemma2}, we have
\begin{align}\label{0627-lll}
e\sigma_{m}(x)=&\frac{e/2}{1+x}+\displaystyle\sum^{m}_{k=2}\int^{1}_{0}\frac{g(s)}{s^{2}}\left(\frac{s}{1+x}\right)^{k}{\rm d}s\nonumber\\
=&\frac{e}{2(1+x)}+\int^{1}_{0}\frac{g(s)}{s^{2}}\displaystyle\sum^{m}_{k=2}\left(\frac{s}{1+x}\right)^{k}{\rm d}s\\
=&\frac{e}{2(1+x)}+\int^{1}_{0}\frac{g(s)}{(1+x)(1+x-s)}{\rm d}s-\int^{1}_{0}\frac{g(s)}{s(1+x-s)}\left(\frac{s}{1+x}\right)^{m}{\rm d}s\nonumber
\end{align}
and
\begin{align}\label{0627-llll}
eS_{m}(x)=&e\displaystyle\sum^{m}_{k=1}\frac{d_{k}}{(x+\frac{11}{12})^{k}}=\frac{e}{2(x+\frac{11}{12})}+e\displaystyle\sum^{m}_{k=2}\frac{d_{k}}{(x+\frac{11}{12})^{k}}\nonumber\\
=&\frac{e}{2(x+\frac{11}{12})}+\int^{1}_{0}\frac{g(t)}{t}\displaystyle\sum^{m}_{k=2}\frac{12}{(12x+11)^k}(12t-1)^{k-1}{\rm d}t+\displaystyle\sum^{m}_{k=2}\frac{\frac{1}{12^{k-1}}(-1)^{k-1}}{(x+\frac{11}{12})^{k}}\\
=&\frac{e}{2(x+\frac{11}{12})}+\int^{1}_{0}\frac{g(t)(t-\frac{1}{12})}{t(1+x-t)(x+\frac{11}{12})}
\left(1-\left(\frac{t-\frac{1}{12}}{x+\frac{11}{12}}\right)^{m-1}\right){\rm d}t\nonumber\\
&+\displaystyle\sum^{m}_{k=2}(-1)^{k-1}\frac{12}{(12x+11)^{k}}.\nonumber
\end{align}
Subtracting \eqref{0627-llll} from \eqref{0627-lll} yields
\begin{align}\label{L-Y-1}
eS_{m}(x)-e\sigma_{m}(x)=&\frac{e}{2}\left(\frac{1}{x+\frac{11}{12}}-\frac{1}{x+1}\right)+\int^{1}_{0}\frac{g(t)(t-\frac{1}{12})}{t(1+x-t)(x+\frac{11}{12})}{\rm d}t
-\int^{1}_{0}\frac{g(t)}{(1+x-t)(x+1)}{\rm d}t\nonumber\\
&+\int^{1}_{0}\frac{g(t)}{t(1+x-t)}\left(\frac{t}{1+x}\right)^{m}{\rm d}t-\int^{1}_{0}\frac{g(t)}{t(1+x-t)}\left(\frac{t-\frac{1}{12}}{x+\frac{11}{12}}\right)^{m}{\rm d}t\nonumber\\
&+\displaystyle\sum^{m}_{k=2}(-1)^{k-1}\frac{12}{(12x+11)^{k}}.
\end{align}

(i)\ If $m\geq1$ is odd. Since $\int^{1}_{0}\frac{g(s)}{s}{\rm d}s=\frac{e}{2}-1$, it follows from
\eqref{L-Y-1} that
\begin{align*}
eS_{m}-e\sigma_{m}(x)=&\left(\frac{1}{x+\frac{11}{12}}-\frac{1}{x+1}\right)\int^{1}_{0}\frac{g(t)}{t}{\rm d}t+\frac{1}{x+\frac{1}{12}}-\frac{1}{x+1}\\
                     &+\int^{1}_{0}\frac{g(t)(t-\frac{1}{12})}{t(1+x-t)(x+\frac{11}{12})}{\rm d}t-\int^{1}_{0}\frac{g(t)}{(1+x)(1+x-t)}{\rm d}t\\
                     &+\int^{1}_{0}\frac{g(t)}{t(1+x-t)}\left(\frac{t}{1+x}\right)^{m}{\rm d}t-\int^{1}_{0}\frac{g(t)}{t(1+x-t)}\left(\frac{t-\frac{11}{12}}{x+\frac{11}{12}}\right)^{m}{\rm d}t\\
                     &+\displaystyle\sum^{m}_{k=2}(-1)^{k-1}\frac{12}{(12x+11)^{k}}.\\
=&\frac{1}{x+\frac{11}{12}}-\frac{1}{x+1}+\displaystyle\sum^{m}_{k=2}(-1)^{k-1}\frac{12}{(12x+11)^{k}}\\
                     &+\int^{1}_{0}\frac{g(t)}{t(1+x-t)}\left(\left(\frac{t}{1+x}\right)^{m}-
                     \left(\frac{t-\frac{1}{12}}{x+\frac{11}{12}}\right)^{m}\right){\rm d}t\\
                     >&\int^{1}_{0}\frac{g(t)}{t(1+x-t)}\left(\left(\frac{t} {1+x}\right)^{m}-\left(\frac{t-\frac{1}{12}}{x+\frac{11}{12}}\right)^{m}\right){\rm d}t\\
                     =&\int^{1}_{0}\frac{g(t)}{t(1+x-t)}\xi^{m-1}\left(\frac{t}{x+1}-\frac{t-\frac{1}{12}}{x+\frac{11}{12}}\right){\rm d}t\\
                     =&\int^{1}_{0}\frac{g(t)\xi^{m-1}}{t(x+1)(12x+11)}{\rm d}t\\
                     >&0.
\end{align*}
Hence, we have $eS_{m}-e\sigma_{m}(t)>0$ for all $m\geq1$ is odd and all $x>0$.

(ii) Suppose $m \geq 2$ is even. By Theorem \ref{444}, the inequality $S_m(t) > \sigma_m(t)$ holds for all $m \geq 6$ and $t > 0$. Therefore, to establish the result for all even $m \geq 2$, it suffices to verify the cases $m = 2$ and $m = 4$. Recall the coefficient values:
\begin{align}\label{eq:coeffs}
b_1 = \frac{1}{2},\
b_2 = \frac{1}{24},\
b_3 = \frac{1}{48},\
b_4 = \frac{73}{5670}
\end{align}
and
\begin{align}\label{eq:coeffs-1}
d_1 = \frac{1}{2},\
d_2 = 0,\
d_3 = \frac{5}{288},\
d_4 = \frac{139}{17280}.
\end{align}
We now verify the inequalities for $x > 0$.

\noindent\textbf{Case $m=2$:}
From \eqref{LLLL-1} with $m=2$, \eqref{eq:coeffs} and \eqref{eq:coeffs-1}, we have 
\begin{align*}
S_2(x) - \sigma_2(x) 
= \frac{d_1}{x + \frac{11}{12}} - \frac{b_1}{x + 1} - \frac{b_2}{(x + 1)^2}=\frac{1/2}{x + \frac{11}{12}} - \frac{1/2}{x + 1} - \frac{1/24}{(x + 1)^2}.
\end{align*}
Combining terms over a common denominator yields
\begin{equation}\label{eq:m2_result}
S_2(x) - \sigma_2(x) = \frac{1}{288(1 + x)^2\left(x + \frac{11}{12}\right)}.
\end{equation}
Since the denominator is strictly positive for $x > 0$, we conclude
$$
S_2(x) - \sigma_2(x) > 0 \ \text{for all} \ x > 0.
$$

\noindent\textbf{Case $m=4$:}
From \eqref{LLLL-1} with  $m=4$, \eqref{eq:coeffs} and \eqref{eq:coeffs-1}, we have 
\begin{align*}
S_4(x) - \sigma_4(x) 
&= \frac{d_1}{x + \frac{11}{12}} + \frac{d_3}{\left(x + \frac{11}{12}\right)^3} + \frac{d_4}{\left(x + \frac{11}{12}\right)^4} \\
&\quad - \frac{b_1}{x + 1} - \frac{b_2}{(x + 1)^2} - \frac{b_3}{(x + 1)^3} - \frac{b_4}{(x + 1)^4} \\
&= \frac{1/2}{x + \frac{11}{12}} + \frac{5/288}{\left(x + \frac{11}{12}\right)^3} + \frac{139/17280}{\left(x + \frac{11}{12}\right)^4} \\
&\quad - \frac{1/2}{x + 1} - \frac{1/24}{(x + 1)^2} - \frac{1/48}{(x + 1)^3} - \frac{73/5670}{(x + 1)^4}.
\end{align*}
Now we prove that for all $x > 0$, the following inequality holds:
\[
S_4(x) - \sigma_4(x) > 0.
\]

{\bf Step 1: Variable Substitution.}
Let $t = x + \frac{11}{12}$. Then for $x > 0$, we have $t > \frac{11}{12}$, and $x + 1 = t + \frac{1}{12}$. 

Rewriting the expression in terms of $t$:
\[
f(t) = \frac{1}{2t} + \frac{5}{288t^3} + \frac{139}{17280t^4} - \frac{1}{2(t+\frac{1}{12})} - \frac{1}{24(t+\frac{1}{12})^2} - \frac{1}{48(t+\frac{1}{12})^3} - \frac{73}{5670(t+\frac{1}{12})^4}.
\]

{\bf Step 2: Common Denominator.}
Express all terms with common denominator $t^4(t+\frac{1}{12})^4$:
\[
f(t) = \frac{N(t)}{t^4(t+\frac{1}{12})^4}
\]
where the numerator polynomial is:
\[
\begin{aligned}
N(t) = &\frac{1}{2}t^3(t+\frac{1}{12})^4 + \frac{5}{288}t(t+\frac{1}{12})^4 + \frac{139}{17280}(t+\frac{1}{12})^4 \\
&-\frac{1}{2}t^4(t+\frac{1}{12})^3 - \frac{1}{24}t^4(t+\frac{1}{12})^2 - \frac{1}{48}t^4(t+\frac{1}{12}) - \frac{73}{5670}t^4.
\end{aligned}
\]

{\bf Step 3: Polynomial Expansion.}
Expand each component:

1. Expand $(t+\frac{1}{12})^n$:
\[
\begin{aligned}
(t+\frac{1}{12})^4 &= t^4 + \frac{1}{3}t^3 + \frac{1}{24}t^2 + \frac{1}{432}t + \frac{1}{20736}; \\
(t+\frac{1}{12})^3 &= t^3 + \frac{1}{4}t^2 + \frac{1}{48}t + \frac{1}{1728}; \\
(t+\frac{1}{12})^2 &= t^2 + \frac{1}{6}t + \frac{1}{144}.
\end{aligned}
\]

2. Compute each term's contribution:
\begin{align*}
\begin{aligned}
&\frac{1}{2}t^3(t+\frac{1}{12})^4 = \frac{1}{2}t^7 + \frac{1}{6}t^6 + \frac{1}{48}t^5 + \frac{1}{864}t^4 + \frac{1}{41472}t^3; \\
&\frac{5}{288}t(t+\frac{1}{12})^4 = \frac{5}{288}t^5 + \frac{5}{864}t^4 + \frac{5}{6912}t^3 + \frac{5}{124416}t^2 + \frac{5}{5971968}t; \\
&\frac{139}{17280}(t+\frac{1}{12})^4 = \frac{139}{17280}t^4 + \frac{139}{51840}t^3 + \frac{139}{414720}t^2 + \frac{139}{7464960}t + \frac{139}{358318080}; \\
&-\frac{1}{2}t^4(t+\frac{1}{12})^3 = -\frac{1}{2}t^7 - \frac{1}{8}t^6 - \frac{1}{96}t^5 - \frac{1}{3456}t^4; \\
&-\frac{1}{24}t^4(t+\frac{1}{12})^2 = -\frac{1}{24}t^6 - \frac{1}{144}t^5 - \frac{1}{3456}t^4; \\
&-\frac{1}{48}t^4(t+\frac{1}{12}) = -\frac{1}{48}t^5 - \frac{1}{576}t^4; \\
&-\frac{73}{5670}t^4.
\end{aligned}
\end{align*}

{\bf Step 4: Term Consolidation.}
Combine like terms:
\begin{align*}
\begin{aligned}
t^7 &: \frac{1}{2} - \frac{1}{2} = 0; \\
t^6 &: \frac{1}{6} - \frac{1}{8} - \frac{1}{24} = 0; \\
t^5 &: \frac{1}{48} + \frac{5}{288} - \frac{1}{96} - \frac{1}{144} - \frac{1}{48} = 0; \\
t^4 &: \frac{1}{864} + \frac{5}{864} + \frac{139}{17280} - \frac{1}{3456} - \frac{1}{3456} - \frac{1}{576} - \frac{73}{5670}\\
&= \frac{2100}{1814400}+\frac{10500}{1814400}+\frac{14630}{1814400}-\frac{525}{1814400}-\frac{525}{1814400}-\frac{3150}{1814400}-\frac{23360}{1814400} \\
&= \frac{2100 + 10500 + 14630 - 525 - 525 - 3150 - 23360}{1814400} \\
&= \frac{2670}{1814400} = \frac{89}{60480} > 0;\\
t^3 &: \frac{1}{41472} + \frac{5}{6912} + \frac{139}{51840}>0;\\
t^2&: \frac{5}{124416} + \frac{139}{414720}>0;\\
t&:\frac{5}{5971968} + \frac{139}{7464960}>0;\\
\text{Constant Term}&: \frac{139}{358318080}.
\end{aligned}
\end{align*}

{\bf Step 5: Conclusion.}
Since:
\begin{itemize}
\item The denominator $t^4(t+\frac{1}{12})^4 > 0$ for $t > \frac{11}{12}.$
\item The numerator $N(t) > 0$ for $t > \frac{11}{12}$ (as all non-zero coefficients are positive).
\end{itemize}

Therefore, for all $x > 0$:
\[
S_4(x) - \sigma_4(x) = f(t) = \frac{N(t)}{t^4(t+\frac{1}{12})^4} > 0.
\]
\\ \\

\textbf{Author contributions} Yongbing Luo wrote the main manuscript text and  Yan Pin provided analysis methods. All authors reviewed the manuscript..

\textbf{Funding} This work was supported by  Supported by Huzhou Key Laboratory of Data Modeling and Analysis (Grant 2023-04).

\textbf{Data Availability} Data sharing not applicable to this article as no datasets were generated or analysed during the current study.

\textbf{Declarations}

\textbf{Conflict of interest}    The authors declare no conflict of interest.

\textbf{Compliance with Ethical Standards} The research does not involve Human Participants or Animals. All authors have approved the contents of this paper and have agreed to the Potential Analysis's submission policies.






\Addresses



\begin{thebibliography}{99}



\bibitem{L-1}
C. Morfici, Y. Hu, On some conversenve to the constant e and
improvements of Carleman's inequality. Carpath, J. Math. 31,
249-254 (2015).
\bibitem{L-5}
X. Yang, On Carleman's inequality. J. Math. And. Appl. 253,
691-694 (2001).
\bibitem{L-3}
M. Gyllenberg, P. Yan, On a conjecture by Yang. J. Math. And.
Appl. 264, 687-690 (2001).



\bibitem{L-4}
Y. Hu, C. Morfici, On the coefficients of an expansion of
$\left(1+\frac{1}{x}\right)^x$ related to Carleman's inequality (
2014). arXiv: 1401.2236 [math.CA].




\bibitem{L-6}
B. Ren, X. Li. Some inequalities related to two expansions
of $\left(1+\frac{1}{x}\right)^x$, 393 (2015).




\end{thebibliography}
\end{document}